\documentclass[12pt,a4paper]{amsart}
\usepackage{graphicx}
\usepackage{ifpdf}
\ifpdf
\fi

\newcommand{\bC}{\mathbb{C}}
\newcommand{\bQ}{\mathbb{Q}}
\newcommand{\bZ}{\mathbb{Z}}
\newcommand{\fg}{\mathfrak{g}}
\newcommand{\fso}{\mathfrak{so}}
\newcommand{\fsp}{\mathfrak{sp}}
\newcommand{\fsl}{\mathfrak{sl}}
\newcommand{\End}{\mathrm{End}}
\newcommand{\Hom}{\mathrm{Hom}}

\newtheorem{prop}{Proposition}

\begin{document}
\title{$R$-matrices for the adjoint representations of $U_q(\fso(n))$}
\author{Deepak Parashar}
\address{%
Mathematics Institute\\
University of warwick\\
Coventry CV4 7AL}
\email{D.Parashar@warwick.ac.uk}
\author{Bruce W. Westbury}
\address{%
Mathematics Institute\\
University of warwick\\
Coventry CV4 7AL}
\email{Bruce.Westbury@warwick.ac.uk}
\date{June 2009}
\begin{abstract} We calculate a two-parameter $R$-matrix which specialises
to the trigonometric $R$-matrix of the minimal affinisation of the adjoint
representation of each of the classical simple Lie algebras
$\fso(n)$ and $\fsp(n)$.
\end{abstract}
\maketitle

\section{Introduction}
The aim of this paper is to present a calculation of a sequence of $R$-matrices.
First we explain our conventions on $R$-matrices. Let $K$ be a field and $V$
a finite dimensional vector space over $K$. Then an $R$-matrix is an element
$R(x)$ of the $K$-algebra $\End(V\otimes V)(x)$. This means it can be represented
by a matrix whose entries are elements of $K(x)$. Put $R_1(x)=R(x)\otimes 1_V$
and $R_2(x)=1_V\otimes R(x)$. Then $R(x)$ is a rational $R$-matrix if it satisfies
the Yang-Baxter equation
\begin{equation*}
R_1(x)R_2(x+y)R_1(y)=R_2(y)R_1(x+y)R_2(x)
\end{equation*}
Similarly, an element $R(u)$ is a trigonometric $R$-matrix if it satisfies
the Yang-Baxter equation
\begin{equation*}
R_1(u)R_2(uv)R_1(v)=R_2(v)R_1(uv)R_2(u)
\end{equation*}

The construction of $R$-matrices is closely related to the representation theory
of quantum groups. In particular, for every finite dimensional complex simple
Lie algebra, $\fg$, there is a Yangian $Y(\fg)$ and also the quantum affine
algebra $U_q(\widehat{\fg})$.

A Yangian $Y$ is a Hopf algebra over $\bQ$ and for each finite dimensional
representation $V$ of $Y$ there is a rational $R$-matrix $R(x)\in
\End(V\otimes V)(x)$. Also, an affine quantised enveloping algebra $U$ is a
Hopf algebra over $\bQ(q)$ and for each finite dimensional representation $V$ of $U$
there is a trigonometric $R$-matrix $R(u)\in \End(V\otimes V)(u)$. These are related
since $R(x)$ can be obtained from $R(u)$ by first substituting $u=q^x$ and then taking
the limit $q\rightarrow 1$.

This is a general construction but very few of these $R$-matrices are known explicitly.
If $\fg$ is a simple Lie algebra then there is an inclusion of Hopf algebras
$U(\fg)\hookrightarrow Y(g)$ of the enveloping algebra in the Yangian and an inclusion
of Hopf algebras $U_q(\fg)\hookrightarrow U_q(\widehat{\fg})$. Thus a finite dimensional
representation $V$ of $Y(\fg)$ or $U_q(\widehat{\fg})$ decomposes as a direct sum
of highest weight representations. In the simplest case the restriction of $V$ is a
highest weight representation. In this case, the rational and trigonometric $R$-matrices
are known explicitly from the tensor product graph method. The only other examples of
explicit $R$-matrices that are known are the rational $R$-matrices of \cite{MR1103907}; the
trigonometric $R$-matrix in \cite{MR1477384} and \cite{MR2123822}; and the trigonometric
$R$-matrix in \cite{MR1123615}.

The representations associated to the rational $R$-matrices in \cite{MR1103907} are the direct
sum of the adjoint representation and the trivial representation for each simple Lie
algebra. In this paper we give the trigonometric $R$-matrices for these representations
for the simple Lie algebras $\fso(n)$ and $\fsp(n)$. These representations of the quantum affine
algebra $U_q(\widehat{\fg})$ are given explicitly in \cite{MR2272099}. This representation is
the minimal affinisation of the adjoint representation (except for $\fsl(n)$).

This calculation also applies in other situations in which the centraliser algebras are quotients of the Birman-Wenzl algebras. Examples are the orthosymplectic super Lie algebras and the twisted affine Lie algebras $A_n^{(2)}$ in \cite{MR806529} and \cite{MR824090}. The twisted affine Lie algebras $A_n^{(2)}$ are also discussed in \cite{MR1401199}.

The general method for finding trigonometric $R$-matrices is to solve the Jimbo equations.
These equations can be solved in special cases using the tensor product graph method
described in \cite{MR1103067}, \cite{MR1126645}, \cite{MR1912580}, \cite{MR1306881}.
This method applies to representations of an affine Kac-Moody algebra whose restriction
to a simple Lie subalgebra is irreducible. Given such a representation, $V$, the vertices
of the graph are the composition factors of $V\otimes V$. The graph is bipartite and the
decomposition of the vertices corresponds to the decomposition of $V\otimes V$ into the
symmetric and anti-symmetric tensors. In the untwisted case, there is an edge connecting
vertices $U_1$, $U_2$ if $\Hom(U_1\otimes U_2, \fg)\ne 0$ where $\fg$ is the adjoint representation.
This method does not apply to the representations we consider as these restrict to the sum of an
irreducible representation and the trivial representation. Nevertheless this method does seem
to give partial information. More specifically, if we follow the tensor product graph method
and construct a graph whose vertices are the composition factors of $V\otimes V$ which have
multiplicity one then the results agree with the calculations in this paper.

\section{$R$-matrices and fusion}
The $R$-matrices in this paper are calculated using the fusion procedure described
in \cite{MR1017340}. This is similar to the cabling procedure in knot theory and the bootstrap
method in quantum field theory. This method was used in \cite{MR1154034}. 
Usually the fusion procedure is applied when the $R$-matrix at a particular value
of the spectral parameter is a projection onto an irreducible representation.

Let $S(k)$ be the group of permutations of $1,2,\ldots ,k$ and for $1\le i\le k-1$
let $s_i$ be the permutation $(i,i+1)$. These are the standard generators. Any word
in these generators can be drawn as a string diagram and we say that a word is reduced
if this string diagram has the property that any two strings cross at most once.
Any permutation can be represented by a reduced word and usualy in many ways. A
fundamental result known as Matsumoto's theorem says that two reduced words represent
the same permutation if and only if they are related by  a finite sequence of moves
of the form
\begin{equation*}
s_is_j \leftrightarrow s_js_i\qquad s_is_{i+1}s_i\leftrightarrow s_{i+1}s_is_{i+1}
\end{equation*}
where $1\le i,j\le k-1$ and $|i-j|>1$.

Next we assume we are given a trigonometric $R$-matrix $R(u)\in \End(V\otimes V)(u)$.
Then we construct a linear operator on $\otimes^k V$ depending on formal parameters
$u_1,\ldots ,u_n$ for each element of $S(k)$. Take a reduced word and the corresponding
string diagram. Then each crossing is labelled by a generator $s_i$. Also each string
is labelled by taking the position of the starting point of the string. If a crossing
is labelled by $s_i$ and the two strings are labelled by $r$ and $s$ then we label the
crossing by $R_i(u_s/u_r)$. Now compose these operators in the order in which they 
appear in the reduced word. This gives a linear operator for each reduced word. It follows from
Matsumoto's theorem and the Yang-Baxter equation that this operator depends only on the
permutation and not on the choice of reduced word.

Now choose four formal parameters which satisfy $u_2/u_1=u_4/u_3=v$ and $u_3/u_1=u_4/u_2=u$
and for fixed $v$ write $S(u)$ for the operator on $\otimes^4V$ associated to the
permutation $(4,3,2,1)$. Put $W=V\otimes V$ and fix $v$. Then $S(u)\in\End(W\otimes W)(u)$
is a trigonometric $R$-matrix. In order to check that the Yang-Baxter equation is satisfied,
it is sufficient to observe that both sides of the equation give an operator which is the
operator associated to a reduced word representing the permutation $(5,6,3,4,1,2)$.

Furthermore, if $R(v)$ is diagonalisable but not invertible then we can write $V\otimes V$
as $M\oplus J$ where $J\ne 0$ is the kernel of $R(v)$ and $M$ is the image. Then we can 
consider $S(u)$ as an element of $\End(M\otimes M)(u)$ and this is also a trigonometric
$R$-matrix.

The notation we use is
\begin{equation*}
\left[an+bx+c\right]=\frac{Q^au^bq^c - Q^{-a}u^{-b}q^{-c}}{q-q^{-1}}
\end{equation*}
This extends the usual definition of quantum integers and becomes equivalent if
we put $Q=q^n$ and $u=q^x$. We will also use the notation
\begin{equation*}
\left\{an+bx+c\right\}={Q^au^bq^c + Q^{-a}u^{-b}q^{-c}}
=\frac{[2an+2bx+2c]}{[an+bx+c]}
\end{equation*}

\section{Orthogonal groups}
Fix $n$ and consider the simple Lie algebra $\fso(n)$. Then irreducible representations
are parametrised by their highest weight and we label the fundamental weights so that
$V(0)$ is the trivial representation, $V(\omega_1)$ is the vector representation and
$V(\omega_2)$ is the adjoint representation. Then
\begin{equation*}
V(\omega_1)\otimes V(\omega_1)\equiv V(0)\oplus V(\omega_2)\oplus V(2\omega_1)
\end{equation*}
The trigonometric $R$-matrix for $V(\omega_1)$ is given explicitly in \cite{MR824090}, \cite{MR1064744}
and is also given by the tensor product graph method. Furthermore, for a certain value
of $v$, $R(v)\in\End(V(\omega_1)\otimes V(\omega_1))(u)$ has image $V(\omega_2)\oplus V(0)$
and kernel $V(2\omega_1)$. Therefore the fusion procedure can be used to find the
trigonometric $R$-matrix for $V(\omega_2)\oplus V(0)$.

This fusion procedure as described is in the centraliser algebra. Instead of using the
centraliser algebras directly we use the Birman-Wenzl-Murakami algebras, $A(k)$,
introduced independently in \cite{MR992598} and \cite{MR1079902}. The advantage of this method is
that we then have a single calculation involving $n$ as a formal parameter instead of a
sequence of calculations depending on the positive integer $n$.

Let $B(k)$ be the Artin braid group with standard generators $\sigma_i$ for
$1\le i\le k-1$. Then the algebra $A(k)$ is a quotient of the group algebra
of $B(k)$ over the field $\bQ(q,Q)$. Define $u_i$ by
\begin{equation*}
u_i=1- \left(\frac{\sigma_i-\sigma_i^{-1}}{q-q^{-1}}\right)
\end{equation*}
Then these satisfy the following tangle relations
\begin{align*}
 \sigma_i\sigma_{i+1}\sigma_i&=\sigma_{i+1}\sigma_i\sigma_{i+1} \\
u_iu_{i\pm 1}\sigma^{\pm 1}&= u_i\sigma_{i\pm 1}^{\mp 1} \\
\sigma^{\pm 1}u_{i\pm 1}u_i&= \sigma_{i\pm 1}^{\mp 1}u_i \\
u_i\sigma_{i\pm 1}\sigma_i&= u_iu_{i\pm 1}\\
\sigma_i\sigma_{i\pm 1}u_i&=u_{i\pm 1}u_i
\end{align*}
and the relations
\begin{align*}
u_i^2 &= \left(1+\frac{Qq^{-1}-Q^{-1}q}{q-q^{-1}}\right)u_i \\
u_i\sigma_{i\pm 1}u_i &= Q^{\pm 1}q^{\mp 1} u_i
\end{align*}
The trigonometric $R$-matrix is given by
\begin{equation*}
R_i(u) = (u-1)q^{-1}\sigma_i - (1-u^{-1})Q^{-1}q\sigma_i^{-1}-(q-q^{-1})(Q^{-1}q-q^{-1})
\end{equation*}
This $R$-matrix satisfies unitarity and crossing symmetry. This means
\begin{align*}
R_i(u)R_i(u^{-1})&=Q[n/2+x-1][n/2-x-1][1+x][1-x] \\
u_iR_{i+1}(u)&=u_iu_{i+1}R_i(Q^{-1}q^2u^{-1})
\end{align*}
Also $R(q^{-2})$ has image $V(\omega_2)\oplus V(0)$ and kernel $V(2\omega_1)$.
Applying the fusion procedure gives the trigonometric $R$-matrix
\begin{equation*}
R_{2i-1}(q^{-2})R_{2i+1}(q^{-2})R_{2i}(q^{-2}u)R_{2i-1}(u)R_{2i+1}(u)R_{2i}(q^2u)
\end{equation*}
This satisfies the Yang-Baxter equation but does not satisfy unitarity and so needs
to be modified. Define $E_i$ by
\begin{equation*}
 (q+q^{-1})E_i = q-s_i+(qQ^{-1}+q^{-1})\frac{u_i}{[n-1]+1}
\end{equation*}
and then we define the trigonometric $R$-matrix by
\begin{equation*}
S_i(u^{1/2}) = 
\frac1k E_{2i-1}E_{2i+1}R_{2i}(q^{-2}u)R_{2i-1}(u)R_{2i+1}(u)R_{2i}(q^2u)
\end{equation*}
where $k=[x][x-1][x+n/2][x+n/2-1]^2$.
Then since $E_i$ can be written as a polynomial in $R_i(q^{-2})$ and $E_i$ is idempotent
it follows that this satisfies the Yang-Baxter equation and unitarity.

The following table shows the various notations for the irreducible representations of $A(4)$.
\begin{center}
\begin{tabular}{c|cccccccc}
& $\emptyset$ & $2$ & $1,1$ & $4$ & $3,1$ & $2,2$ & $2,1,1$ & $1^4$ \\
$\fso(n)$ & $0$ & $\omega_2$ & $2\omega_1$ & $\omega_4$ & $\omega_3+\omega_1$ & $2\omega_2$ & $\omega_2+2\omega_1$ & $4\omega_1$ \\
$\fsp(n)$ & $0$ & $2\omega_1$ & $\omega_2$ & $4\omega_1$ & $\omega_2+2\omega_1$ & $2\omega_2$ & $\omega_3+\omega_1$ & $\omega_4$ \\
\end{tabular}
\end{center}

Then the eigenvalues for the rank one idempotents are as follows:
\begin{center}
\begin{tabular}{c|l}
$1^4$ & $\quad [x-1][x-2][n/2+x-2]$ \\
$2,1,1$ & $-[x-1][x+2][n/2+x-2]$ \\
$2,2$ & $\quad [x+1][x+2][n/2+x-2]$ \\
$2$ & $\quad [x-1][x+2][-n/2+x+2]$ \\
\end{tabular}
\end{center}

\begin{prop}
The $2\times 2$ matrix associated to the empty partition is $A=(a_{ij})$
where
\begin{multline*}
a_{11}=\{2x\}\{n/2-1\}[2][n/2-2]\frac{[n+x-2]}{[n+x/2-1]} \\
+[x+2][x-1][x+n/2-2]
\end{multline*}
\begin{multline*}
a_{22}=\{-2x\}\{n/2-1\}[2][n/2-2]\frac{[n-x-2]}{[n-x/2-1]} \\
-[x-2][x+1][x-n/2+2]
\end{multline*}

\begin{equation*}
a_{12}= [2x][2]\frac{[n/2-2]^2}{[x+n/2-1]}
\end{equation*}

\begin{equation*}
a_{21}=\frac{[2x][n/2][n-1]\{n/2-2\}}{\{n/2-1\}^2[n/2+x-1]}
\end{equation*}
\end{prop}

The $3\times 3$ matrix associated to (the adjoint representation) is
$B=(b_{ij})$ where
\begin{align*}
b_{11}=b_{22}&=[2][n/2-2]\frac{\{n/2+x-1\}}{\{n/2-1\}} \\
b_{13}=b_{23}&=\{n/2\}\{n/2-3\}[n/2-1][x] \\
b_{31}=b_{32}&=[2]^2[x][n/2-2] \\
b_{12}=b_{21}&
\end{align*}

In order to compare with \cite[\S 5.2 (21)]{MR1103907} we replace $B$ by $P^{-1}BP$ where 
\[ P=\begin{pmatrix}
1 & 0 & 1 \\ 1 & 0 & -1 \\ 0 & 1 & 0 
     \end{pmatrix} \]

\begin{prop} The matrix $P^{-1}BP$ is
\[ \begin{pmatrix}
c_{11} & c_{21} & 0 \\ c_{12} & c_{22} & 0\\ 0 & 0 & c_{33} 
     \end{pmatrix} \]
where 
\begin{align*} c_{11}&= [x+2][x-1][x+n/2-2]+2[2][n/2-2]\frac{\{x+n/2-1\}}{\{n/2-1\}} \\
c_{22}&= [-x+2][-x-1][-x+n/2-2]+2[2][n/2-2]\frac{\{-x+n/2-1\}}{\{n/2-1\}} 
\end{align*}
\begin{align*}
 c_{21}&= [2]^2[x][n/2-2] \\
 c_{12}&= 2[x][n/2-1]\{n/2\}\{n/2-3\} \\
 c_{33}&= -[x+2][x-1][x+n/2-2]
\end{align*}
Note the relation $c_{11}(u)=c_{22}(u^{-1})$ is satisfied.
\end{prop}

\section{Twisted $R$-matrices}

Let $L$ be a simple Lie algebra with a diagram automorphism $\sigma$ of order $m$.
This gives a grading of $L$, $L=\oplus_{i=0}^{m-1}$, by $\bZ/m\bZ$. In particular
$L_0$ is the subalgebra fixed by $\sigma$. Then we have three affine Kac-Moody algebras
\[ \widehat{L}^{(1)}, \widehat{L}^{(k)}, \widehat{L_0}^{(1)} \]
and two inclusions
\begin{equation}\label{homs}
 \widehat{L}^{(k)}\rightarrow\widehat{L}^{(1)},\qquad \widehat{L}^{(1)}_0\rightarrow \widehat{L}^{(1)}
\end{equation}
Each of these subalgebras is the subalgebra fixed by an automorphism of order $k$
and the restriction of each automorphism to $L$ is $\sigma$.

Associated to each affine Kac-Moody algebra, $\widehat{L}^{(k)}$, is the Drinfeld-Jimbo quantum group,
$U_q(\widehat{L}^{(k)})$. Let $V$ be a finite dimensional irreducible representation of
$U_q(\widehat{L}^{(k)})$. Then there is a trigonometric $R$-matrix with spectral parameter acting
on $V\otimes V$ which can be found by solving the Jimbo equations.

The inclusions in \eqref{homs} have quantum analogues
\begin{equation}
 U_q(\widehat{L}^{(k)})\rightarrow U_q(\widehat{L}^{(1)}),
\qquad U_q(\widehat{L}^{(1)}_0)\rightarrow U_q(\widehat{L}^{(1)})
\end{equation}

Let $U$ and $V$ be finite dimensional representations of $U_q(\widehat{L}^{(1)})$
whose restrictions to $U_q(\widehat{L}^{(k)})$ and $U_q(\widehat{L}^{(1)}_0)$
are irreducible. Then we have three trigonometric $R$-matrices which we denote by
\[ R^{(1)}_{UV}, S^{(k)}_{UV}, S^{(1)}_{UV} \]

The simplest case is when $U$ and $V$ restrict to an irreducible representation of $L_0$.
In this case all three $R$-matrices can be found by the tensor product graph method.
This case is studied in \cite{MR1401199}.

The case we are interested is when the restrictions of $U$ and $V$ to representations
of $L_0$ are $L_1\oplus \bC$. In these cases the restrictions of $U$ and $V$ to representations
of $L$ are irreducible and so the $R$-matrix $R^{(1)}_{UV}$ can be found by the tensor product
graph method.

There are two examples where we can apply our calculations.
\subsection{$A_{2k}^{(2)}$}
In this example we have $L=A_{2k}=\fsl(2k+1)$ and $L_0=B_k=\fso(2k+1)$.
The vector representation of $\fsl(2k+1)$ (and its dual) restricts to the vector representation
of $\fso(2k+1)$. This is one of the examples studied in \cite{MR1401199}.

The representation $V(2\omega_1)$ of $A_{2k}$ is the symmetric square of the vector representation.
This representation (and its dual) restricts to the representation $V(2\omega_1)\oplus\bC$ of $B_k$.

The tensor product graph for $V(2\omega_1)\otimes V(2\omega_1)$ is
\[ \includegraphics{ad.1} \]

The tensor product graph for $V(2\omega_1)\otimes V(2\omega_{2k})$ is
\[ \includegraphics{ad.3} \]

There is no tensor product graph for the other two $R$-matrices. However for the twisted $R$-matrix
we expect to see
\[ \includegraphics{ad.5} \]
and for the untwisted $R$-matrix we expect to see
\[ \includegraphics{ad.7} \]

The branching rules in this case are given in \cite[\S 11.9 Theorem II]{MR0002127}.
\subsection{$A_{2k-1}^{(2)}$}
In this example we have $L=A_{2k-1}=\fsl(2k)$ and $L_0=C_k=\fsp(2k)$.

The vector representation of $\fsl(2k)$ (and its dual) restricts to the vector representation
of $\fsp(2k)$. This is one of the examples studied in \cite{MR1401199}.

The representation $V(\omega_2)$ of $A_{2k-1}$ is the exterior square of the vector representation.
This representation (and its dual) restricts to the representation $V(\omega_2)\oplus\bC$ of $C_k$.

The tensor product graph for $V(\omega_2)\otimes V(\omega_2)$ is
\[ \includegraphics{ad.2} \]

The tensor product graph for $V(\omega_2)\otimes V(\omega_{2k-2})$ is
\[ \includegraphics{ad.4} \]

There is no tensor product graph for the other two $R$-matrices. However for the twisted $R$-matrix
we expect to see
\[ \includegraphics{ad.6} \]
and for the untwisted $R$-matrix we expect to see
\[ \includegraphics{ad.8} \]

\section{Implications}

The motivation for this work comes from \cite{MR1960703} where it was observed that the rational
$R$-matrices in \cite{MR1103907} can be written in terms of co-ordinates arising from the
universal Lie algebra.

Then the rank one idempotents are given by
\begin{center}
\begin{tabular}{c|l}
$1^4$ & $\quad [x+\alpha][x-\beta][x+\gamma]$ \\
$2,1,1$ & $-[x+\alpha][x+\beta][x+\gamma]$ \\
$2,2$ & $\quad [x-\alpha][x-\beta][x+\gamma]$ \\
$2$ & $\quad [x+\alpha][x-\beta][x-\gamma]$ \\
\end{tabular}
\end{center}

The $2\times 2$ matrix associated to the empty partition is $A=(a_{ij})$
where
\begin{multline*}
a_{11}=-\{2x\}\{\alpha+\beta+\gamma\}[\alpha][\beta][\gamma]\frac{[n+x-2]}{[n+x/2-1]} \\
+[x+\alpha][x+\beta][x+\gamma]
\end{multline*}
\begin{multline*}
a_{22}=\{2x\}\{\alpha+\beta+\gamma\}[\alpha][\beta][\gamma]\frac{[n-x-2]}{[n-x/2-1]} \\
+[-x+\alpha][-x+\beta][-x+\gamma]
\end{multline*}

\begin{equation*}
a_{12}= [2x][2]\frac{[n/2-2]^2}{[x+n/2-1]}
\end{equation*}

\begin{equation*}
a_{21}=\frac{[2x][n/2][n-1]\{n/2-2\}}{\{n/2-1\}^2[n/2+x-1]}
\end{equation*}

For the adjoint representation the $3\times 3$ matrix is
\[ \begin{pmatrix}
c_{11} & c_{21} & 0 \\ c_{12} & c_{22} & 0\\ 0 & 0 & c_{33} 
     \end{pmatrix} \]
where 
\begin{align*} c_{11}&= [x+\alpha][x+\beta][x+\gamma]-2[\alpha][\beta][\gamma]\frac{\{x+\alpha+\beta+\gamma\}}{\{\alpha+\beta+\gamma\}} \\
c_{22}&= [-x+\alpha][-x+\beta][-x+\gamma]-2[\alpha][\beta][\gamma]\frac{\{-x+\alpha+\beta+\gamma\}}{\{\alpha+\beta+\gamma\}} 
\end{align*}

\begin{align*}
 c_{21}&= [2]^2[x][\alpha+\beta+\gamma-1] \\
 c_{12}&= 2[x][\alpha+\beta+\gamma]\{\alpha+\beta\}\{\alpha+\gamma\}\{\beta+\gamma\} \\
 c_{33}&= -[x+\alpha][x+\beta][x+\gamma]
\end{align*}
Note the relation $c_{11}(u)=c_{22}(u^{-1})$ is satisfied.

\bibliographystyle{halpha}
\bibliography{adjoint}
\end{document}